\documentclass [12pt,a4paper,reqno]{amsart}
\usepackage{amsmath,amsthm}
\usepackage{amsfonts}
\usepackage{amssymb}
\allowdisplaybreaks

\newcommand{\be}{\begin{equation}}
\newcommand{\ef}{\end{equation}}
\chardef\bslash=`\\ 

\newtheorem*{thm*}{Theorem}

\theoremstyle{definition}

\newtheorem*{remark*}{Remarks}
\newtheorem*{defn*}{Definition}

\theoremstyle{remark}

\newcommand{\G}{\Gamma}
\newcommand{\wt}{\widetilde}
\newcommand{\wh}{\widehat}
\newcommand{\fc}{\frac}

\newcommand{\iy}{\infty}
 \renewcommand{\sectionmark}[1]{}

\newcommand{\Be}{Beltrami}

\newcommand{\qc} {quasiconformal}

\newcommand{\ve}{\varepsilon}

\newcommand{\Te} {Teichm\"{u}ller}

 \usepackage{amsfonts}
\newcommand{\field}[1]{\mathbb{#1}}
\newcommand{\g}{\gamma}

\newcommand{\D}{\field{D}}
\newcommand{\om}{\omega}

\newcommand{\ov}{\overline}
\newcommand{\vp}{\varphi}
\newcommand{\hC}{\wh{\field{C}}}
\newcommand{\C}{\field{C}}

\newcommand{\B}{\mathbf{B}}
\newcommand{\T}{\mathbf{T}}
\newcommand{\Belt}{\operatorname{Belt}}

\newcommand{\Teich}{\operatorname{Teich}}
\newcommand{\dist}{\operatorname{dist}}

\renewcommand{\a} {\alpha}

\newcommand{\ld}{\lambda}
\newcommand{\kp}{\kappa}

\begin{document}

\title{A general coefficient theorem for univalent functions}
\author{Samuel L. Krushkal}

\begin{abstract} Using the Bers isomorphism theorem for Teichm\"{u}ller spaces of punctured \linebreak 
Riemann surfaces and some of their other complex geometric features, we prove a general theorem on maximization of 
homogeneous polynomial (in fact, more general
holomorphic) coefficient functionals $J(f) = J(a_{m_1},
a_{m_2},\dots, a_{m_n}) $ on some classes of univalent functions in the unit disk naturally connected with the canonical class $S$. 
The given functional $J$ is lifted to the Teichm\"{u}ller space $\mathbf T_1$ of the punctured disk $\D_{*} = \{0 < |z| < 1\}$ which is biholomorphically equivalent to the Bers fiber space over the universal Teichm\"{u}ller space. This  generates a  positive subharmonic function on the disk $\{|t| < 4\}$ with 
$\sup_{|t|<4} u(t) = \max_{\mathbf T_1} |J|$ attaining this  
maximal value only on the boundary circle, which correspond to rotations of the Koebe function.  

This theorem  implies new sharp distortion estimates for univalent
functions giving explicitly the extremal functions, and creates a new bridge between Teichm\"{u}ller space theory and geometric
complex analysis. In particular, it provides an alternate and direct proof of the Bieberbach conjecture.

\end{abstract}

\date{\today\hskip4mm({GenCoefThm.tex})}

\maketitle

\bigskip

{\small {\textbf {2010 Mathematics Subject Classification:} Primary: 30C50, 30C75, 30F60; 
Secondary 30C55, 30C62, 31A05, 32L05, 32Q45}

\medskip

}

\bigskip

\markboth{Samuel L. Krushkal}{A general coefficient theorem}
\pagestyle{headings}

\bigskip\bigskip
\centerline{\bf 1. PREAMBLE. CLASSES OF FUNCTIONS}

\bigskip\noindent
{\bf 1.1}.  Estimating holomorphic functionals on the classes of
univalent functions depending on the Taylor coefficients of these
functions is important in various geometric and physical
applications of complex analysis, because these coefficients reflect
the fundamental intrinsic features of conformal maps. One of the
main collections of univalent functions is the canonical class $S$ formed by
univalent functions on the unit disk $\D = \{|z| < 1\}$ with $f(0) = 0, \ f^\prime(0) = 1$.

We consider a family of univalent functions in $\D$ closely related to the class $S$. This family, denoted by $\wh S(1)$, is the completion
in the topology of locally uniform convergence on $\D$ of the set of
univalent functions $f(z) = a_1 z + a_2 z^2 + \dots$ with $|a_1| =
1$, having quasiconformal extensions across the unit circle $\mathbb S^1 =
\partial \D$ to the whole sphere $\hC = \C \cup \{\iy\}$ which
satisfy $f(1) = 1$.

Equivalently, this family is a disjunct union
$$
\wh S(1) = \bigcup_{- \pi \le \theta < \pi} S_\theta(1),
$$
where $S_\theta(1)$ consists of univalent functions
$$
f(z) = e^{i \theta} z + a_2 z^2 + \dots
$$
with quasiconformal extensions to $\hC$ satisfying $f(1) = 1$ (also
completed in the indicated weak topology).

Every $f \in S$ has its representative $\wh f$ in $\wh S(1)$ 
(not necessarily unique) 
obtained by pre and post compositions of $f$ with rotations $z
\mapsto e^{i \alpha} z$ about the origin. 
In the general case, the equality $f(1) = 1$ holds in the sense of the boundary correspondence under conformal maps, i.e., in terms of the Carath\'{e}odory prime ends. 

It suffices to verify this for $f \in S$ mapping $\D$ onto a Jordan
domain $G$. Since $f(0) = 0, \ |f^\prime(0)| = 1$, it follows from
the Schwarz lemma that $L = \partial G$ must intersect $\mathbb S^1$. Take
one of the intersection points $e^{i \theta} = f(z_0) \in L\cap
\mathbb S^1$; then the function
  \be\label{1}
f_{\tau, \theta}(z) =  e^{- i \theta} f(e^{i \tau} z) \quad
\text{with} \ \ \tau = \arg z_0
\end{equation}
belongs to $\wh S(1)$. The converse is trivial: given a function 
$f \in \wh S(1)$ belonging to some $S_\theta(1)$, then 
$f(e^{- i \theta} z) \in S$; for other $f \in S$, this follows by weak approximation. 

This relation implies, in particular, that the functions conformal
in the closed disk $\ov \D$ are dense in each class $S_\theta(1)$.
Such a dense subset is formed, for example, by the images of the
homotopy functions $[f]_r(z) = \fc{1}{r} f(r z)$ with real $r \in
(0, 1)$.

The inversions $F(z) = 1/f(1/z)$ of $f \in \wh S(1)$ form the
corresponding classes $\Sigma_\theta(1)$ of nonvanishing univalent
functions on the complementary disk $\D^* = \{z \in \hC: \ |z|
> 1\}$, with expansions
$$
F(z) =  e^{- i \theta} z + b_0 + b_1 z^{-1} + b_2 z^{-2} + \dots,
\quad  F(1) = 1,
$$
and $\wh \Sigma(1) = \bigcup_\theta \Sigma_\theta(1)$. 

We use the standard notation $\Sigma$ for the class of all 
univalent functions on $\D^*$ with expansions 
$$
F(z) = z + b_0 + b_1 z^{-1} + b_2 z^{-2} + \dots \ . 
$$

The condition $w(1) = 1$ is added to determine uniquely the map $w$
by its Schwarzian
$$
S_w(z) = (w^{\prime\prime}(z)/w^\prime(z))^\prime -
(w^{\prime\prime}(z)/w^\prime(z))^2/2, \quad z \in \D,
$$
as well as by its Beltrami coefficient $\mu_w(z) = \partial_{\ov z}
w/\partial_z w$ on $\C$ and to ensure the holomorphic dependence of
$w$ from $S_w$ and $\mu_w$. 

Simple computations yield that the coefficients $a_n$ of $f \in
S_\theta(1)$ and the corresponding coefficients $b_j$ of $F(z) = 1/f(1/z) \in \Sigma_\theta(1)$ are related by
$$
b_0 + e^{2i \theta} a_2 = 0, \quad b_n + \sum \limits_{j=1}^{n}
\epsilon_{n,j}  b_{n-j} a_{j+1} + \epsilon_{n+2,0} a_{n+2} = 0,
\quad n = 1, 2, ... \ ,
$$
where $\epsilon_{n,j}$ are the entire powers of $e^{i \theta}$. This successively implies the representations of $a_n$ by $b_j$ via
 \be\label{2}
a_n = (- 1)^{n-1} \epsilon_{n-1,0}  b_0^{n-1} - (- 1)^{n-1} (n - 2)
\epsilon_{1,n-3} b_1 b_0^{n-3} + \text{lower terms with respect to}
\ b_0
\end{equation}
(so $a_2 = - e^{- 2i \theta} b_0, \  a_3 = - e^{- 2i \theta} b_1 +
e^{- 2i \theta} b_0^2, \ \dots$).

\bigskip\noindent 
{\bf 1.3}. Our goal is to estimate on these
classes the general nonconstant homogeneous polynomial coefficient functionals
 \be\label{3}
J(f) = J(a_{m_1}, a_{m_2},\dots, a_{m_s}), \quad J(\mathbf 0) = 0,
\end{equation}
where $2 \le m_1 < m_2 < \dots < m_s \le n$ and $\mathbf 0 = (0,
\dots, 0)$ is the origin of $\C^s$. Such $J$ generate by (2) the corresponding polynomial functionals
$$
\wt J(F) = \wt J(b_0, \dots , b_{m_{n-2}})
$$
on $\wh \Sigma(1)$. The existence of extremal functions of $J(f)$ and $\wt J(F)$ follows from compactness of these classes in the topology of locally uniform convergence.

It is well known that many distortion functionals on $S$
attain their maximal value on the Koebe function
 \be\label{4}
\kp_0(z) = \fc{z}{(1 - z)^2} = z + \sum\limits_2^\iy n z^n
\end{equation}
mapping the unit disk onto the complement of the ray $\{w = t : \
-1/4 \le t \le \iy\}$, and on rotations $k_\theta(z) = e^{- i
\theta} \kp_0(e^{i \theta} z)$ of this function.

Our main theorem sheds light on this phenomenon and shows the  
crucial role of two distinguished sets intrinsically connected with functionals on $\wh S(1)$: the zero set 
$$ 
\mathcal Z_J = \{f \in \wh S(1): J(f) = 0\}
$$ 
(which is a subset of the pluripolar set of plurisubharmonic function $\log |J(a_{m_1}, a_{m_2},\dots, a_{m_s})|$ on $\C^s$) 
and the set of rotations
 \be\label{5}
\mathcal K = \{\kp_{\tau, \theta}(z) = e^{- i \theta} \kp_0(e^{i
\tau} z)\}.
\end{equation}

\bigskip\bigskip
\centerline{\bf 2. MAIN THEOREM AND ITS CONSEQUENCES}

\bigskip\noindent
{\bf 2.1. Main Theorem}. The main basic result of this paper is

\bigskip\noindent {\bf Theorem 1}. {\it Any homogeneous polynomial functional (3), whose zero set $\mathcal Z_J$ is separated from the set (5), is maximized on $\wh S(1)$ only by the functions $f_0 \in
\mathcal K$. }

\bigskip
In other words, {\it any extremal function  $f_0$ of any homogeneous coefficient
functional $J$ on the class $\wh S(1)$ is  the Koebe function (4)
composed with pre and post rotations about the origin, unless}
$J(f_0) = 0$. The assumption $\mathcal Z_J \cap \mathcal K =
\emptyset$ cannot be omitted.

\bigskip
The proof of this theorem involves a deep result from Teichm\"{u}ller space theory given by the Bers isomorphism theorem \cite{Be}. The functional $J$ is lifted from $S_\theta(1)$ to the Teichm\"{u}ller space $\T_1$ of the punctured disk $\D_{*} = \{0 < |z| < 1\}$. This space is biholomorphically equivalent to the Bers
fiber space $\mathcal F(\T)$ over the universal Teichm\"{u}ller
space $\T = \Teich (\D)$. This generates a holomorphic functional $\mathcal J(\vp, t)$ on $\mathcal F(\T)$ covering $J$, with the same range domain. Here $\vp = S_F \in \T$ are the Schwarzian derivatives of functions $F \in \Sigma_\theta(1)$, while the variable $t$ runs
over the fiber domain $F_\vp(\D)$ defined by $\vp$ .

A crucial step is to maximize $|\mathcal J(\vp, t)|$ over $\vp$ by a fixed $t$. We apply an approximation of the underlying space $\T$ by the finitely dimensional Teichm\"{u}ller spaces of the punctured spheres in the weak topology of locally uniform convergence in $\C$. 
Any such space is foliated by Teichm\"{u}ller-Kobayashi geodesic disks. We deal with restrictions of  $\mathcal J(\vp, t)$ to these disks, taking their appropriate dense contable collection. 
This implies a maximal logarithmically subharmonic function $u_\theta(t) > 0$ on a domain located in the disk $\D_4 = \{|t| < 4\}$.

Repeating this construction for all $\theta$, one creates a 
logarithmically subharmonic function $u(t) = \sup u_\theta (t)$ on the disk $\D_4$ with
$$
\max_t u(t) = \sup_{\T_1} |\mathcal J| = \max_S |J|.
$$
This maximal value is attained on the boundary of $\D_4$ whose
points correspond to the function $\kp_0(z)$ composed with rotations.

One can see from the proof that Theorem 1 can be straightforwardly
extended to more general appropriate holomorphic distortion
functionals $J(f)$ on $\wh S(1)$. In particular, it is valid for the linear combinations of polynomial functionals 
$$
\a_1 J_1 + \dots + \a_n J_n
$$ 
with positive $\a_1, \dots, \a_n$. Moreover, one can allow  the terms $\a_j J_j$ of this combination to depend from pairwise different collections of coefficients $a_{m_s}$. 

\bigskip\noindent
{\bf 2.2. Applications}. Theorem 1 implies new sharp distortion
estimates for univalent functions giving explicitly the extremal
functions, and creates a new bridge between geometric complex
analysis and Teichm\"{u}ller space theory.

We illustrate this by two important results. First of all, as a
direct consequence of Theorem 1, we have

\bigskip\noindent
{\bf Corollary 1}. {\it If a functional $J(f)$ on $S$ satisfies
$$
\max_S |J(f)| = \max_{\wh S(1)} |J(f)| \quad\text{and} \ \
|J(\kp_0)| > 0,
$$
then every extremal of this functional is a a rotated Koebe function $\kp_{\theta}$
with some $\theta \in [-\pi, \pi)$. For a homogeneous $J$, any
function $\kp_\theta$ is extremal for $J$. }

\bigskip
The assumptions of Corollary 1 are satisfied, for example, by
monomial functionals
$$
J(f) = a_{m_1}^{q_1} \dots \ a_{m_n}^{q_n},  \quad q_j \in \field N,
$$
and by functionals of the form $J(f) = P(a_n)$, where $P(z)$ is a
polynomial without free term and with positive coefficients,
$$
P(z) = c_1 z + c_2 z^2 + \dots + c_N z^N, \quad c_j > 0 \ \ (n > 1).
$$

As is well known, there were several classical conjectures about the coefficients that have been intensively investigated for a very long time. They include
the Bieberbach conjecture which states that in the class $S$ the coefficients are
estimated by $|a_n| \le n$ with equality only for $\kp_\theta$, as
well as several other well-known conjectures that imply the
Bieberbach conjecture. Most of them have been proved by the de Branges theorem \cite {DB} (see also \cite{Ha}). 

Corollary 1, applied to  $J(f) = a_n$, implies that for all $f \in S$,
$$
\max_S |a_n| = |a_n(\kp_\theta)| = n.
$$
This is an alternate and direct proof of the Bieberbach conjecture.

\bigskip
The second application concerns sharp estimation of the coefficients of
Schwarzian derivatives
$$
S_f(z) = \sum\limits_0^\iy \a_n z^n \quad (|z| < 1) \ \ \text{of}
$$
on $S$. This problem also has been investigated by many authors.

For $f_{\tau, \theta}(z) = e^{- i \theta} f(e^{i \tau} z), \ f \in
S$, we have $S_{f_{\tau, \theta}}(z) = e^{2i \theta} S_f(e^{i \tau}
z)$. Noting also that the Schwarzian of the Koebe function is given
by
$$
S_{\kp_0}(z) = - \fc{6}{(1 - z^2)^2} = - 6 \ \sum\limits_0^\iy (n + 1) z^{2n},
$$
(and is even), one immediately obtains from Corollary 1 that the
even coefficients $\a_{2 n}$ of all Schwarzians $S_f$ on $S$ are
sharply estimated by
 \be\label{6}
|\a_{2 n}(S_f)| \le |\a_{2 n}(\kp_0)| = 6(n + 1),
\end{equation}
with equality only for $f = \kp_\theta$. Together with the Moebius
invariance of $S_f$ and its higher derivatives, the estimate (6) provides the bound of distortion at an arbitrary point of $\D$ (cf. \cite{Ah}, \cite{Ber}, \cite{H}, \cite{Sc}, \cite{Ta}).

\bigskip\bigskip
\centerline{\bf 3. DIGRESSION TO TEICHM\"{U}LLER SPACES}

\bigskip
We briefly recall some needed results from Teichm\"{u}ller space theory on spaces involved in order to prove Theorem 1; the details can be found, for example, in \cite{Be}, \cite{GL}, \cite{Le}.

\bigskip\noindent
{\bf 3.1}.  The universal Teichm\"{u}ller space $\T = \Teich (\D)$ is the space of
quasisymmetric homeomorphisms of the unit circle $\mathbb S^1$ factorized by
M\"{o}bius maps;  all Teichm\"{u}ller spaces have their isometric
copies in $\T$.

The canonical complex Banach structure on $\T$ is defined by
factorization of the ball of the Beltrami coefficients (or complex
dilatations)
$$
\Belt(\D)_1 = \{\mu \in L_\iy(\C): \ \mu|\D^* = 0, \ \|\mu\| < 1\},
$$
letting $\mu_1, \mu_2 \in \Belt(\D)_1$ be equivalent if the
corresponding \qc \ maps $w^{\mu_1}, w^{\mu_2}$ (solutions to the
Beltrami equation $\partial_{\ov{z}} w = \mu \partial_z w$ with $\mu
= \mu_1, \mu_2$) coincide on the unit circle $\mathbb S^1 = \partial \D^*$
(hence, on $\ov{\D^*}$). Such $\mu$ and the corresponding maps
$w^\mu$ are called $\T$-{\it equivalent}. The equivalence classes
$[w^\mu]_\T$ are in one-to-one correspondence with the Schwarzian derivatives
$$
S_w(z) = \left(\frac{w^{\prime\prime}(z)}{w^\prime(z)}\right)^\prime - \frac{1}{2} \left(\frac{w^{\prime\prime}(z)}{w^\prime(z)}\right)^2 \quad (w = w^\mu(z), \ \ z \in \D^*).
$$

Note that for each locally univalent function $w(z)$ on a simply
connected hyperbolic domain $D \subset \hC$, its Schwarzian
derivative belongs to the complex Banach space $\B(D)$ of
hyperbolically bounded holomorphic functions on $D$ with the norm
$$
\|\vp\|_\B = \sup_D \ld_D^{-2}(z) |\vp(z)|,
$$
where $\ld_D(z) |dz|$ is the hyperbolic metric on $D$ of Gaussian
curvature $- 4$; hence $\vp(z) = O(z^{-4})$ as $z \to \iy$ if $\iy
\in D$. In particular, for the unit disk, 
$$
\ld_\D(z) = 1/(1 - |z|^2).
$$

The space $\B(D)$ is dual to the Bergman space $A_1(D)$, a subspace
of $L_1(D)$ formed by integrable holomorphic functions (quadratic differentials $\vp(z) dz^2$ on $D$), since every linear functional $l(\vp)$ on $A_1(D)$ is represented in the form
 \be\label{7}
l(\vp) = \langle \psi, \vp \rangle_D = \iint\limits_D \
\ld_D^{-2}(z) \ov{\psi(z)} \vp(z) dx dy
\end{equation}
with a uniquely determined $\psi \in \B(D)$.

The Schwarzians $S_{w^\mu}(z)$ with $\mu \in \Belt(\D)_1$ range over a bounded domain in the space $\B = \B(\D^*)$. This domain models the space $\T$. It lies in the ball $\{\|\vp\|_\B < 6\}$ and contains the ball $\{\|\vp\|_\B < 2\}$. In this model, the
Teichm\"{u}ller spaces of all hyperbolic Riemann surfaces are
contained in $\T$ as its complex submanifolds.

The factorizing projection
$$
\phi_\T(\mu) = S_{w^\mu}: \ \Belt(\D)_1 \to \T
$$
is a holomorphic map from $L_\iy(\D)$ to $\B$. This map is a split
submersion, which means that $\phi_\T$ has local holomorphic
sections (see, e.g., [GL]).

Note that both equations $S_w = \vp$ and $\partial_{\ov z} w = \mu
\partial_z w$ (on $\D^*$ and $\D$, respectively) determine their
solutions in $\Sigma_\theta$ uniquely, so the values $w^\mu(z_0)$
for any fixed $z_0 \in \C$ and the Taylor coefficients $b_1, b_2,
\dots$ of $w^\mu \in \Sigma_\theta$ depend holomorphically on $\mu
\in \Belt(\D)_1$ and on $S_{w^\mu} \in \T$.

\bigskip\noindent
{\bf 3.2}. The points of Teichm\"{u}ller space $\T_1 =
\Teich(\D_{*})$ of the punctured disk $\D_{*} = \D \setminus \{0\}$
are the classes $[\mu]_{\T_1}$ of $\T_1$-{\it equivalent} \Be \
coefficients $\mu \in \Belt(\D)_1$ so that the corresponding \qc \
automorphisms $w^\mu$ of the unit disk coincide on both boundary
components (unit circle $\mathbb S^1 = \{|z| =1\}$ and the puncture $z = 0$)
and are homotopic on $\D \setminus \{0\}$. This space can be endowed
with a canonical complex structure of a complex Banach manifold and
embedded into $\T$ using uniformization.

Namely, the disk $\D_{*}$ is conformally equivalent to the factor
$\D/\G$, where $\G$ is a cyclic parabolic Fuchsian group acting
discontinuously on $\D$ and $\D^*$. The functions $\mu \in
L_\iy(\D)$ are lifted to $\D$ as the \Be \ $(-1, 1)$-measurable
forms  $\wt \mu d\ov{z}/dz$ in $\D$ with respect to $\G$, i.e., via
$(\wt \mu \circ \g) \ov{\g^\prime}/\g^\prime = \wt \mu, \ \g \in
\G$, forming the Banach space $L_\iy(\D, \G)$.

We extend these $\wt \mu$ by zero to $\D^*$ and consider the unit
ball $\Belt(\D, \G)_1$ of $L_\iy(\D, \G)$. Then the corresponding
Schwarzians $S_{w^{\wt \mu}|\D^*}$ belong to $\T$. Moreover, $\T_1$
is canonically isomorphic to the subspace $\T(\G) = \T \cap \B(\G)$,
where $\B(\G)$ consists of elements $\vp \in \B$ satisfying $(\vp
\circ \g) (\g^\prime)^2 = \vp$ in $\D^*$ for all $\g \in \G$.

Due to the Bers isomorphism theorem, the space $\T_1$ is
biholomorphically isomorphic to the Bers fiber space
$$
\mathcal F(\T) = \{(\phi_\T(\mu), z) \in \T \times \C: \ \mu \in
\Belt(\D)_1, \ z \in w^\mu(\D)\}
$$
over the universal space $\T$ with holomorphic projection $\pi(\psi,
z) = \psi$ (see \cite{Be}).

This fiber space is a bounded hyperbolic domain in $\B \times \C$
and represents the collection of domains $D_\mu = w^\mu(\D)$ as a
holomorphic family over the space $\T$. For every $z \in \D$,  its
orbit $w^\mu(z)$ in $\T_1$ is a holomorphic curve over $\T$.

The indicated isomorphism between $\T_1$ and $\mathcal F(\T)$ is
induced by the inclusion map \linebreak 
$j: \ \D_{*} \hookrightarrow \D$
forgetting the puncture at the origin via
 \be\label{8}
\mu \mapsto (S_{w^{\mu_1}}, w^{\mu_1}(0)) \quad \text{with} \ \
\mu_1 = j_{*} \mu := (\mu \circ j_0) \ov{j_0^\prime}/j_0^\prime, 
\end{equation} 
where $j_0$ is the lift of $j$ to $\D$. 

In the line with our goals, we slightly modified the Bers  construction, applying quasiconformal maps $F^\mu$ of $\D_{*}$ admitting conformal extension to $\D^*$ (and accordingly using  the Beltrami coefficients $\mu$ supported in the disk) (cf. \cite{Kr2}). These changes are not essential and do not affect the underlying features of the Bers isomorphism (giving the same space up to a biholomorphic isomorphism).  

The Bers theorem is valid for Teichm\"{u}ller spaces $\T(X_0 \setminus
\{x_0\})$ of all punctured hyperbolic Riemann surfaces $X_0
\setminus \{x_0\}$ and implies that $\T(X_0 \setminus \{x_0\})$ is
biholomorphically isomorphic to the Bers fiber space $ \mathcal
F(\T(X_0))$ over $\T(X_0)$.

\bigskip\noindent
{\bf 3.3}. We shall consider also the space $\T(0, n)$ of punctured
spheres (Riemann surfaces of genus zero)
$$
X_{\mathbf z} = \hC \setminus \{0, 1, z_1 \dots, z_{n-3}, \iy\}
$$
defined by ordered $n$-tuples $\mathbf z = (0, 1, z_1, \dots,
z_{n-3}, \iy), \ n > 4$ with distinct $z_j \in \C \setminus \{0,
1\}$.

Fix a collection $\mathbf z^0 = (0, 1, z_1^0, \dots, z_{n-3}^0,
\iy)$ with $ z_j^0 \in S^1$ defining the base point $X_{\mathbf
z^0}$ of the Teichm\"{u}ller space $\T(0, n) = \T(X_{\mathbf z^0})$. Its
points are the equivalence classes $[\mu]$ of Beltrami coefficients
from the ball $\Belt(\C)_1 = \{\mu \in L_\iy(\C): \ \|\mu\|_\iy <
1\}$ under the relation: $\mu_1 \sim \mu_2$, if the corresponding
\qc \ homeomorphisms $w^{\mu_1}, w^{\mu_2}: \ X_{\mathbf a^0} \to
X_{\mathbf a}$  are homotopic on $X_{\mathbf a^0}$ (and hence
coincide in the points $0, 1, z_1^0, \dots, z_{n-3}^0, \iy$). This
models $\T(0, n)$ as the quotient space $\T(0, n) = \Belt(\C)_1/\sim
$ with complex Banach structure of dimension $n - 3$ inherited from
the ball $\Belt(\C)_1$.

Another canonical model of $\T(0, n) = \T(X_{\mathbf z^0})$ is
obtained again using the uniformization. The surface $X_{\mathbf
z^0}$ is conformally equivalent to the quotient space $U/\G_0$, where
$\G_0$ is a torsion free Fuchsian group of the first kind acting
discontinuously on $\D \cup \D^*$. The functions $\mu \in
L_\iy(X_{\mathbf z^0})$ are lifted to $\D$ as the Beltrami $(-1,
1)$-measurable forms  $\wt \mu d\ov{z}/dz$ in $\D$ with respect to
$\G_0$ which satisfy $(\wt \mu \circ \g) \ov{\g^\prime}/\g^\prime =
\wt \mu, \ \g \in \G_0$ and form the Banach space $L_\iy(\D, \G_0)$.

After extending these $\wt \mu$ by zero to $\D^*$, the  Schwarzians
$S_{w^{\wt \mu}|\D^*}$ for $\|\wt \mu\|_\iy$  belong to $\T$ and
form its subspace regarded as the {\it Teichm\"{u}ller space
$\T(\G_0)$ of the group $\G_0$}. It is canonically isomorphic to the
space $\T(X_{\mathbf z^0})$, and moreover,
 \be\label{9}
\T(\G_0) = \T \cap \B(\G_0),
\end{equation}
where $\B(\G_0)$ is an $(n - 3)$-dimensional subspace of $\B$ which
consists of elements $\vp \in \B$ satisfying $(\vp \circ \g)
(\g^\prime)^2 = \vp$ for all $\g \in \G_0$ (holomorphic
$\G_0$-automorphic forms of degree $- 4$); see, e.g. \cite{Le}.

This leads to the representation of the space $\T(X_{\mathbf z^0})$
as a bounded domain in the complex Euclidean space $\C^{n-3}$.

Note that $\B(\G_0)$ has the same elements as the space $A_1(\D^*,
\G_0)$ of integrable holomorphic forms of degree $- 4$ with norm
$\|\vp\|_{A_1(\D^*, \G_0)} = \iint_{\D^*/\G_0} |\vp(z)| dx dy$; and
similar to (7), every linear functional $l(\vp)$ on $A_1(\D^*,
\G_0)$ is represented in the form
$$
l(\vp) = \langle \psi, \vp \rangle_{\D/\G_0} :=
\iint\limits_{\D^*/\G_0} (1 - |z|^2)^2 \ \ov{\psi(z)} \vp(z) dx dy
$$
with uniquely determined $\psi \in \B(\G_0)$.

Any Teichm\"{u}ller space is a complete metric space with intrinsic
Teichm\"{u}ller metric defined by quasiconformal maps. By the
Royden-Gardiner theorem, this metric equals the hyperbolic  Kobayashi metric determined by the complex structure (see, e.g., \cite{EKK}, \cite{GL}).

\bigskip\bigskip
\centerline{\bf 4. PROOF OF THEOREM 1}

\bigskip
We accomplish the proof of this theorem in three stages.

\bigskip\noindent
$\mathbf{1^0}$. \ Using the Bers isomorphism theorem, one can regard
the points of the space $\T_1$ as the pairs $X_{F^\mu} = (S_{F^\mu}, F^\mu(0))$, where $\mu \in \Belt(\D)_1$ obey $\T_1$-equivalence (hence, also $\T$-equivalence). Letting
 \be\label{10}
\wh J(\mu) = \wt J(F^\mu),
\end{equation}
we lift the given functional  from the sets $S_\theta(1)$ and
$\Sigma_\theta(1)$ onto the ball $\Belt(\D)_1$. Then, under the
indicated $\T_1$-equivalence, i.e., by the quotient map
$$
\phi_{\T_1}: \ \Belt(\D)_1 \to \T_1, \quad \mu \to [\mu]_{\T_1},
$$
the functional $\wt J(F^\mu)$ is pushed down to a bounded
holomorphic functional $\mathcal J$ on the space $\T_1$ with the same range domain. 

Equivalently, one can apply the quotient map $\Belt(\D)_1 \to \T$ (i.e., $\T$-equivalence) and compose  the descended 
functional on $\T$ with the natural holomorphic map $\iota_1: \ \T_1 \to \T$ generated by the inclusion $\D_{*} \hookrightarrow \D$ forgetting the puncture. Note that since the coefficients $b_0, \ b_1, \dots$ of $F^\mu \in \Sigma_\theta$   are uniquely determined by its Schwarzian $S_{F^\mu}$, the values of $\mathcal J$ in the points $X_1, \ X_2 \in \T_1$ with $\iota_1(X_1) = \iota_1(X_2)$ are equal. 

To avoid a complication of notations, we shall denote the
composition of $\mathcal J$ with biholomorphism $\T_1 \cong \mathcal F(\T)$ again by $\mathcal J$, getting by (8) and (10),
 \be\label{11}
\mathcal J(X_{F^\mu}) = \mathcal J(S_{F^\mu}, \ t), \quad t =
F^\mu(0).
\end{equation}
This yields a plurisubharmonic functional $|\mathcal J(S_{F^\mu}, t)|$ on $\mathcal F(\T)$. 

By Koebe's one-quarter theorem, the boundary of domain
$F(\D^*)$ under any nonvanishing univalent function $F(z)= z + b_0 + b_1 z^{-1} + \dots$ on $\D^*$ (an inversion of $f(z) = z + a_2 z^2 + \dots \in S$) is located in the disk $\{|w - b_0| \le 2\}$ and $|b_0| = |a_2| \le 2$; hence, the variable $t$ in the representation (11) runs over a subdomain $D_\theta$ of the disk $\D_4 = \{|t| < 4\}$.

We define on this domain the function
 \be\label{12}
u_\theta(t) = \sup_{S_{F^\mu}} |\mathcal J(S_{F^\mu}, t)|,
\end{equation}
taking the supremum over all $S_{F^\mu} \in \T$ admissible for a given $t = F^\mu(0) \in D_\theta$, that means over the pairs
$(S_{F^\mu}, t) \in \mathcal F(\T)$ with a fixed $t$.

\bigskip\noindent
$\mathbf{2^0}$. \ Our goal is to establish that this function
inherits subharmonicity, which is given by the following basic
lemma.

\bigskip\noindent
{\bf Lemma 1}. {\it The function $u_\theta(t)$ is subharmonic in the domain $D_\theta$. }

\bigskip
{\bf Proof}. We apply a weak approximation of the underlying space $\T$ (and simultaneously of the space $\T_1$) by finite dimensional Teichm\"{u}ller spaces of the punctured spheres in the topology of locally uniform convergence on $\C$.

Note that the space $\B$ containing $\T$ is not separable;  thus, it does not admit a strong sequential approximation in $\B$-norm, and the weak (locally uniform) convergence of the Schwarzians $S_f$
on $\D$ (together with $f \in \wh S(1)$) also restrains the growth
of $\mathcal J$ on $\mathcal F(\T)$. It ensures, for example, that the maximum of $|\mathcal J|$ cannot increase under semicontinuous regularization in the strong topology.

We take the set of points
$$
E = \{e^{\pi s i/2^n}, \ s = 0, 1, \dots, 2^{n+1} - 1; \ n = 1, 2,
\dots\}
$$
(which is dense on the unit circle) and consider the punctured
spheres
$$
X_m = \hC \setminus \{e^{\pi s i/2^n}, \ s = 0, 1, \dots, 2^{n+1} -
1\}, \quad m = 2^{n+1},
$$
and their universal holomorphic covering maps $g_m: \ \D \to X_m$
normalized by $g_m(0) = 0, \ g_m^\prime(0) > 0$.

The radial slits from the infinite point to all the points $e^{\pi s
i/2^n}$ form a canonical dissection $L_m$ of $X_m$ and define the
simply connected surface $X_m^\prime = X_m \setminus L_m$. Any
covering map $g_m$ determines a Fuchsian group $\G_m$ of covering
transformations acting discontinuosly in both disks $\D$ and $\D^*$.
The (open) fundamental polygons $P_m$ of $\G_m$ in $\D$
corresponding to the dissection $L_m$ is a regular circular
$2^{n+1}$-gon centered at the origin which can be chosen to have a
vertex at $1$. These polygons entirely increase and exhaust the disk
$\D$, and the restriction of $g_m$ to $P_m$ is univalent.

Similarly, we take in the complementary disk $\D^*$ the mirror polygons 
$P_m^*$ and the covering maps $g_m^*(z) = 1/\ov{g_m(1/ \ov z)}$
which define the mirror surfaces $X_m^*$.

Now we approximate the maps $F^\mu \in \Sigma_\theta(1)$ by
homeomorphisms $F^{\mu_m}$ having in $\D = \{|z| < 1\}$ the 
Beltrami coefficients
  \be\label{13}
\mu_m = [g_m]_{*} \mu := (\mu \circ g_m) \ov{g_m^\prime}/g_m^\prime, \ \ n= 1, 2, \dots \ .
\end{equation}
Each $F^{\mu_m}$ is again $k$-quasiconformal  and compatible with
the group $\G_m$. As $m \to \iy$, the coefficients $\mu_m$ are
convergent to $\mu$ almost everywhere on $\C$; thus, the maps $
F^{\mu_m}$ are convergent to $F^\mu$ uniformly in the spherical
metric on $\hC$. Note also that $\mu_m$ depend holomorphically on $\mu$ as elements of $L_\iy$; hence, $F^{\mu_m}(0)$ {\it is a holomorphic function of} $t = F^\mu(0)$. 

In view of importance of this fact, we also provide its 
alternate proof based onlocal existence theorem from \cite{Kr1} which we present here as

\bigskip\noindent
{\bf Lemma 2}. {\it Let $D$ be a finitely connected domain on the Riemann 
sphere $\hC$. Assume that there are a set $E_0$ of positive
two-dimensional Lebesgue measure and a finite number of points
 $z_1, z_2, ..., z_m$ distinguished in $D$. Let
$\a_1, \a_2, ..., \a_m$ be non-negative integers assigned to $z_1,
z_2, ..., z_m$, respectively, so that $\a_j = 0$ if $z_j \in E_0$.

Then, for a sufficiently small $\ve_0 > 0$ and $\varepsilon \in (0,
\varepsilon_0)$, and for any given collection of numbers $w_{sj}, s = 0, 1, ..., \a_j, \ j = 1,2, ..., m$, which satisfy the conditions
$w_{0j} \in D$, \
$$
|w_{0j} - z_j| \le \ve, \ \ |w_{1j} - 1| \le \ve, \ \ |w_{sj}| \le
\ve \ (s = 0, 1, \dots   a_j, \ j = 1, ..., m),
$$
there exists a \qc \ automorphism $h_\ve$ of domain $D$ which is
conformal on $D \setminus E_0$ and satisfies
$$
h_\ve^{(s)}(z_j) = w_{sj} \quad \text{for all} \ s =0, 1, ..., \a_j,
\ j = 1, ..., m.
$$
Moreover, the \Be \ coefficient $\mu_{h_\ve}(z) = \partial_{\bar z}
h_\ve/\partial_z h_\ve$ of $h$ on $E_0$ satisfies $\| \mu_{h_\ve}
\|_\iy \le M \ve$ and depends holomorphically on 
$\{w_{sj}\}$ (with indicated values of $s$ and $ j$).   
The constants $\ve_0$ and $M$ depend only on
the sets $D, E_0$ and the vectors $(z_1, ..., z_m)$ and $(\a_1, ..., \a_m)$. }

\bigskip 
This lemma ensures the existence of local holomorphic sections to the maps $\mu \mapsto F^\mu(z)$ and $\mu \mapsto  
S_{F^\mu}$ defined by maps $\mu  \mapsto 
\{\frac{d^s}{dz^s} F^\mu(z_j)\}$.  

One can fix $\mu_0 \in \Belt(\D)_1$ and apply Lemma 2, for example, to collections 
$$
\{F^{\mu_0}(0), \ F^{\mu_0}(z_{j,n})\}, \quad z_{j,n} = e^{\pi j i/2^n} \in F^{\mu_0}(E),  
$$ 
with prescribed images of these points (sufficiently close to their originals). This provides the corresponding Beltrami coefficients $\mu_h$ depending holomorphically on $t = F^\mu(0)$. Then also the coefficients $\mu_{h \circ f^{\mu_0}}$ are holomorphic in $t$,  and this holomorphy is preserved under transform (13). 

As a result, one obtains that the Beltrami coefficients 
$$
\mu_{h,m} := [g_m]_{*} \mu_h 
$$ 
and the corresponding values  $F^{\mu_{h,m}}(0)$ are holomorphic functions of the variable $t = F^\mu(0)$, 
which completes the proof of the claim.  

\bigskip
We have established that the resulting function $\mathcal J(S_{F^{\mu_m}}, t)$ with $t =
F^{\mu_m}(0)$ is separably holomorphic in its arguments $S_{F^{\mu_m}}$
and $t$; hence, by the Hartogs theorem (extended to Banach spaces) this function is jointly holomorphic in
$(S_{F^{\mu_m}}, t) \in \mathcal F(\T)$. 
  
We now choose in $\T(0, m) \setminus \{\mathbf 0\}$ represented as a subdomain of the space $\B(\G_m)$ a countable dense subset
$$
E^{(m)} = \{\vp_1, \vp_2, \dots, \vp_p, \dots\}. 
$$
For any of its point $\vp_p$,  the corresponding extremal
Teich\"{u}ller disk $\D(\vp_p)$ joining this point with the origin of $\B(\G_m)$ does not meet other points from this set 
(this follows from the uniqueness of Teich\"{u}ller extremal map). Recall also that each disk $\D(\vp_p)$ is formed by the Schwarzians $S_{F^{\tau \mu_{p;m}}}$ with $|\tau| < 1$ and 
$$
\mu_{p;m}(z) = |\psi_{p;m}(z)|/\psi_{p;m}(z)
$$ 
with appropriate $\psi_{p;m} \in A_1(\D, \G_m), \ \|\psi_{p;m}\|_1 = 1$. 

The restrictions of the functional $\mathcal J(S_{F^{\tau \mu_{p;m}}}, t)$ to
these disks are holomorphic functions of $(\tau, t)$; moreover, the above construction provides that all these restrictions are
holomorphic in $t$ in some common domain $D_m \subset \D_4$
containing the point $t = 0$, provided that $|\tau| \le k < 1$. We use the maximal common holomorphy domain; it is located in a disk $\{|t| < r_0\}, \ r_0 \le 4$.  

Maximization over $\tau$ implies the subharmonic functions 
$\log U_{p;m}(t)$ with 
 \be\label{14}
U_{p;m}(t) = \sup_{|\tau| <1} |\mathcal J(S_{F^{\tau \mu_{p;m}}}, t)| \quad (t = F^{\mu_{p;m}}(0), \ \ p = 1, 2, \dots) 
\end{equation}
in the domain $D_m$. We consider the upper envelope of this sequence
 \be\label{15}
u_m(t) = \sup_p U_{p;m}(t) 
\end{equation}
defined in some domain $D_m \subset \D_4$ containing the origin, and take the upper semicontinuous regularization 
of this envelope 
 \be\label{16}
u_m(t) = \limsup\limits_{t^\prime \to t} u_m(t^\prime)
\end{equation}
(by abuse of notation, we shall denote the regularizations by the same letter as the original functions). 

The general properties of subharmonic functions in the Euclidean spaces imply that such a regularization also is subharmonic; therefore, the function (16) is logarithmically subharmonic in the domain $D_m$.

On the other hand, in view of the density of the set $E^{(m)}$ in
$\T(0, m) = \T(\G_m)$ and of compactness of the spaces $\T$ and
$\T_1$ in the topology of locally uniform convergence on $\D^*$, we also have the equality
$$
u_m(t)= \sup \ \{|\mathcal J(S_{F^{\mu_m}}, t)|: \ S_{F^{\mu_m}} \in \T(\G_m)\} 
$$
with $\mu_m$ given by (13) and $t = F^{\mu_m}(0)$. 

Repeating the above construction for each $m = 4, 5, \dots$, one obtains the corresponding increasing sequences of the uniformizing groups $\G_m$, domains $D_m$ arising from fibers and of maximal logarithmically subharmonic functions $u_m$, so that
$$
\G_m \subset \G_{m+1} \subset \G_{m+2} \subset \dots \ , \ \ D_m
\subset D_{m+1} \subset D_{m+2} \subset \dots \subset \D_4,
$$
and for a fixed $m_0 \ge 4$ the functions $u_m$ with $m \ge m_0$
satisfy
$$
u_m(t) \le u_{m+1}(t) \le u_{m+2}(t) \le \dots \
$$
for $t \in D_{m_0}$.

Every function $u_m$ satisfies on $D_{m_0}$ the mean value inequality
$$
u_m(t_0) \le \frac{1}{2 \pi} \int\limits_0^{2 \pi} u_m(t_0 + r e^{i \theta})| d \theta  
$$ 
(for $0 < r < \dist(t_0, \partial D_{m_0})$). Since this sequence is  monotone increasing, one derives that this inequality is preserved also for the limit function 
$$
\wt u_\theta(t) = \lim\limits_{m\to \iy} u_m(t). 
$$
Its regularization via (16) 
 \be\label{17}
u_\theta(t) = \limsup\limits_{t^\prime \to t} \wt u_\theta(t^\prime)  
\end{equation}
provides a positive upper semicontinuous and hence logarithmically subharmonic function on the domain
$$
D_\theta = \bigcup_m D_m \subset \D_4
$$
filled by the corresponding values of $F^\mu(0)$.

Taking into account that the approximating sequence  $F^{[g_m]_{*} \mu}$ constructed above  for every $F^\mu \in \Sigma_\theta(1)$ is convergent to this $F^\mu$ locally uniformly on $\C$, and this construction can be made 
also for the functions $F^\mu$ close to any $1/\kp_{\tau,
\theta}(z)$, we obtain that the boundary of the limit domain
$D_\theta$ has at least one common point with the circle $\{|t| = 4\}$.

On the other hand, both limit operations (16), (17) involve the weakly convergent appropriate maximizing sequences of the Schwarzians $S_{F_j} \in \T$ on $\D^*$, and the weak compactness of the class
$\wh \Sigma_\theta(1)$ implies that the maximal limit function (17) must coincide with the function defined by (12). This  completes the proof of Lemma 1.

\bigskip\noindent
$\mathbf{3^0}$. \ Applying the above constructions to all $\theta, \ -\pi \le \theta < \pi$, and maximizing in $\theta$, one obtains in a similar manner to the above a maximal function
 \be\label{18}
u(t) = \sup_\theta u_\theta(t),
\end{equation}
which is defined and subharmonic in the domain 
$$ 
D = \bigcup_\theta D_\theta \subset \D_4.
$$
We need to establish that this domain must coincide with the disk $\D_4$. This is a consequence of the following facts. 

First, noting that in view of (1) the Schwarzians of two  corresponding functions $F_1 \in \Sigma_{\theta_1}(1)$ and $F_2 \in \Sigma_{\theta_2}(1)$ are related by 
$$ 
S_{F_1}(z) = S_{F_2} \circ \g_{1,2}(z) \g_{1,2}^\prime(z)^2, 
$$ 
and hence both belong to the underlying universal Teichm\"{u}ller space $\T = \Teich (\D)$, similarly as above one can apply the $\T_1$-equivalence on dense subset $\wh \Sigma^0(1)$ in the union $\wh \Sigma(1)$ formed by all 
univalent functions 
$$
F_\theta(z) = e^{-i \theta} z + b_0 + b_1 z^{-2} + \dots
$$ 
on $\D^*$ (preserving $z = 1$) with quasiconformal extension to $\hC$. This implies the (non-contractible) universal quotient space 
$$
\mathcal T = \wh \Sigma^0(1)/\thicksim.
$$ 

It is holomorphically equivalent to the quotient  
$\Belt(\D)_1/\thicksim$ 
under the same equivalence relation. 
The Bers isomorphism $\beta: \T_1 \to \mathcal F(\T)$ defined by (8) is lifted to this space as a holomorphic surjection of $\mathcal T$ onto $\mathcal F(\T)$. 

In a similar way, one can apply the above construction from the previous step simultaneously to any finite union of the quotient spaces 
 \be\label{19}
\mathcal T_m = \bigcup_{j=1}^m \ \wh \Sigma_{\theta_j}^0/\thicksim \ = \bigcup_{j=1}^m \{(S_{F_{\theta_j}}, F_\theta^\mu(0)) \} \ \simeq \T_1 \cup \dots \cup \T_1,    
\end{equation} 
where the equivalence relation $\thicksim$ again means $\T_1$-equivalence and 
$$
\mathbf F_\theta^\mu(0) := (F_{\theta_1}^{\mu_1}(0), \dots , F_{\theta_m}^{\mu_m}(0)). 
$$ 
The Beltrami coefficients  $\mu_j \in \Belt(\D)_1$ are chosen here independently.    
The corresponding collection 
$\beta = (\beta_1, \dots, \beta_m)$  
of the Bers isomorphisms 
$$
\beta_j: \ \{(S_{F_{\theta_j}}, F_{\theta_j}^{\mu_j}(0))\} \to \mathcal F(\T)
$$  
determines a holomorphic surjection of the space $\mathcal T_m$ onto $\mathcal F(\T)$.   

We select a dense countable subset $\{ \theta_j\} \subset [-\pi, \pi]$, and consider the increasing unions (19) for $\theta_j$ from this set (which corresponds to taking the maximum in (18)).   

The spaces $\mathcal T_m$ possess an important property 
inherited from the classes  $\Sigma$ and $\wh \Sigma(1)$ inherited from the class $\Sigma$ and preserved under connection (1). This is a circular symmetry of maps $F_{\tau,\sigma}$, since for any $F \in \Sigma$, all functions $e^{- i \a} F(e^{i \a}z)$ also belong to this class. This generates the complex homotopy 
$$
F_r(z) = r F(z/r): \ \D^* \times \D \to \hC \quad (F \in \Sigma)         
$$ 
and the corresponding homotopy on the universal Teichm\"{u}ller space $\T$. 

To preserve the fixed point $1$, we choose such $r \in \overline{\D}$, for which the maps  $F_{r;\tau, \theta} \in \wh \Sigma^0(1)$ (related by (1) to $F_r \in \Sigma$) satisfy $F_{r;\tau, \theta}(1) = 1$. The Beltrami coefficients of these maps are connected by   
$$ 
\mu_r(z) = \mu(z/r) e^{2 i \theta^\prime} \ov z/z, 
$$ 
where $\mu$ is the Beltrami coefficient of the initial function  $F$ on $\D$ and $\theta^\prime$ depends only on 
$\arg r$ and $\theta$ coming from (1).  
The indicated relation for Beltrami coefficients implies  the desired circular symmetry of the corresponding functions $F_{r;\tau, \theta}$. 

This yields that any union (19) contains together with every of its points $(S_{F_\theta}, F_\theta^\mu(0))$ also the appropriate points 
$(S_{F_{r;\tau, \theta'}^{\mu_r}}, F_{r;\tau, \theta'}^{\mu_r}(0))$ with $|r| = 1$ and  $|F_{r;\tau, \theta'}^{\mu_r}(0)| = |F_\theta^\mu(0)|$ (i.e., with the second coordinate located on a circular arc centred at $F_\theta^\mu(0)$). 

Now the assumption of homogeneity of the original functional $J$ implies for these points the crucial equality 
  \be\label{20}
|\mathcal J(S_{F_{r;\tau, \theta'}^{\mu_r}}, F_{r;\tau, \theta'}^{\mu_r}(0))| = 
|\mathcal J(S_{F^\mu}, F^\mu(0))|        
\end{equation} 
(provided that the points are outside of the polar set of function $\log |\mathcal J|$ in $\T_1$). 
This equality implies that the value domains of $\mathcal J$ on the spaces $\mathcal T_m$ admit some circular symmetry.    

Finally, the assumption of the theorem, that the original functional $J$ does not vanish on the set (5) of functions $\kp_{\tau, \theta}$, provides that the above arguments  
are valid also  for the maps $F^\mu \in \wh \Sigma^0(1)$ approximating the corresponding rotations $F_{\tau,\theta}$  of the function $F_0(z) = 1/\kp_0(1/z)$. 

We combine all this with another important property of the  function (4). Namely, by Koebe's one-quarter theorem mentioned above, any univalent function 
$$
w = F(z) = z + b_0 + b_1 z^{-1} + \dots \in \Sigma
$$ 
maps the disk $\D^*$ onto a domain  $F(\D^*)$ whose boundary is located in the disk $\{|w - b_0| \le 2\}$. In addition, for all $F$ inverting the univalent $f \in S$, we have the equality $b_0 = a_2$. 

This implies that $|b_0| \le 2$ and that the complementary domain to $F(D^*)$ is located in the disk $\{|w| \le 4\}$. 
Also it follows from above that the boundary points of this disk (with $|w| = 4$) are covered only by functions 
$$
F(z) = \fc{1}{\kp_\theta(z)} = z + 2 e^{i \theta} + \fc{e^{2i \theta}}{z}, \quad - \pi \le \theta < \pi.
$$

Accordingly, in the case of the corresponding maps $F_{\tau,\theta} \in \wh \Sigma(1)$, we have that the points of the circle $\{|w| = 4\}$ are filled only by compositions of the function $F(z) = 1/\kp_0(z)$ with pre and post rotations about the origin. 

Hence, letting $m \to \infty$, one obtains in the limit that  the function $u(t)$ given by (18) is defined and subharmonic on the whole disk $\{|t| < 4\}$; moreover, in view of the density of the set $\{\theta_j\}$, the relation (20) yields that this function must be radial. 

On the other hand, the weak compactness of the class $\wh S(1)$ implies that 
$$
\limsup\limits_{|t| \to 4} u(t) \ge \inf_\theta |J(\kp_{\tau,
\theta})| > 0.     
$$
Since $u(0) = J(\mathbf 0) = 0$, the function $u(t)$ is not constant. As a radial function, it must be monotone increasing on $[0, 4]$ and hence attains its maximal value only on the boundary circle $\{|t| = 4\}$, which is equivalent to the assertion of the theorem.

\bigskip\bigskip 
\centerline{\bf 5. ADDITIONAL REMARKS}

\bigskip\noindent
{\bf 1. A modified proof of Lemma 1}. 
In view of the crucial role of the basic Lemma 1, we 
provide a somewhat modified proof combining the applied interpolation of $\T$ and $\T_1$ by finite dimensional Teichm\"{u}ller spaces with a quasiconformal surgery of maps $F$. 
Namely, for each $\mu \in \Belt(\D)_1$ there exists a number $m_0$ depending on $[\mu]_{\T_1}$ such that for all $m \ge m_0$ both points $F^\mu(0)$ and $F^{\mu_m}(0)$ are located in the same polygon $F^\mu(P_m)$. For all $F^\mu(0)$ with $\|\mu\|_\iy \le k < 1$  (which corresponds to $S_{F^\mu}$ running in the hyperbolic ball in $\T$ of radius $\tanh^{-1} k$ centered  at the origin), one can take the same lower bound $m_0 = m_0(k)$. 

We join these points by a hyperbolic segment $l_m$ in $F^\mu(P_m)$ and consider its $\delta$-neighborhood 
$U_\delta = \{z: \ \dist(z, l_m) < \delta\} \Subset F^\mu(P_m)$ with some $\delta > 0$. 

Now, let $\om_m$ be the quasiconformal automorphism of $F^\mu(P_m)$ identical on the boundary of this polygon and moving $F^{\mu_m}(0)$ into $F^\mu(0)$ (the extremal dilatation of such $\om_m$ is given explicitly in \cite{Te}). Then
$$
\om_m \circ F^{\mu_m}|_{\partial P_m} = F^{\mu_m}|_{\partial P_m},
$$
and this equality extends by periodicity to $\D$. As $m \to \infty$, the sequence $\{\om_m \circ F^{\mu_n}\}$ also is convergent to
$F^\mu$ locally uniformly in the spherical metric. This surgery
implies for each point $F^\mu(0)$ a finite dimensional manifold
$M_m$ in $\T_1$ over this point representing $\T(0, m)$(which
depends only on the class $[\mu]_{\T_1}$). 
The resulting function $\mathcal J(S_{F^{\mu_m}}, t)$ with $t =
F^\mu(0)$ is jointly holomorphic in 
$(S_{F^{\mu_m}}, t) \in \mathcal F(\T)$. 

Now, assuming that $|\tau| \le k < 1$, we use instead of (14) the functions 
$$
U_{m,k,p}(t) = \sup_{|\tau| \le k} |\mathcal J(S_{F^{\tau \mu_p}}, t)| \quad (t = F^\mu(0), \ \ p = m_0(k), m_0(k) + 1, \dots), 
$$
which are defined and subharmonic in a common domain $D_m(k) \subset \D_4$. It is located in a disk $\{|t| < r_0(k)\}, \ 
r_0(k) < 4$ and contains the point $t = 0$.  

Similar to (16), we take the upper envelope 
$u_m(t) = \sup_p U_p(t)$ followed by its upper semicontinuous regularization.

\bigskip\noindent
{\bf 2}.  Theorem 1 implies that in fact the maximization of 
homogeneous polynomial coefficient functionals on $\wh S(1)$ is reduced to
maximization of trigonometric polynomials generated by the extremals
$\kp_{\tau, \theta}$ of $J$, which is much more elementary and is solved by techniques of calculus.

\bigskip\noindent
{\bf 3. A distortion theorem for the class $\wh S(1)$}. 
Pick a generic polynomial 
$P(z) = \sum\limits_1^N c_n z^n$ with zero free term and 
positive $c_n$, and consider on $\wh S(1)$ the polynomial functional of the form
  \be\label{21}
J(f) = P \Bigl(\fc{a_n}{n}\Bigr) = \fc{c_1}{n} a_n + \fc{c_2}{n^2} a_n^2 + \dots + \fc{c_N}{n^N} a_n^N \quad (n > 1). 
\end{equation}
By Theorem 1 for all $f \in \wh S(1)$ we have the
estimate $|a_n| \le n$ with equality for $f = \kp_{\tau, \theta}$.
It remains to find the unknown values of $\tau, \ \theta$
determining these extremals.

This implies, for example, the following distortion
theorem for the class $\wh S(1)$.

\bigskip\noindent 
{\bf Theorem 2}. {\it Let the polynomial $P(z)$ 
have on the unit circle exactly one critical point $z_0$. Then 
the correspomding functional (21) is estimated by 
 \be\label{22}
\max_{\wh S(1)} |J(f)| = |J(\kp_{\tau_0, \theta_0})|,
\end{equation}
where $\tau_0$ and $\theta_0$ are defined (in general, not uniquely) from the equations}
 \be\label{23}
\fc{d |P(e^{i(\om + \tau)})|}{d \om} \Big\vert_{\om = \arg z_0} = 0, \quad \kp_{\tau_0, \theta_0}(1) = 1.
\end{equation}

\bigskip\noindent
{\bf Proof}. For any given $J(f)$, Theorem 1 implies the sharp
estimate (22) with some $\tau_0, \ \theta_0 \in [- \pi, \pi]$.
Since, by assumption, $|P(z)|$ has on the unit circle $\mathbb S^1$  exactly one critical point $z = z_0$, we have
 \be\label{24}
\max |P(e^{i\om})| = |P(z_0)| 
\end{equation} 
(note that under this assumption $|P(z)|$ must have on $\mathbb S^1$ a simple zero). 

The equality (24) implies the first equation in (23). This equation determines the admissible values of $\tau_0$; hence, the functional $J$ attains its maximum on the set $\{\kp_{\tau_0,\theta}\}$.
Now, to find its extremal functions $\kp_{\tau_0, \theta_0}$, one must solve the equation $\kp_{\tau_0, \theta}(1) = 1$, 
completing the proof. 

A similar estimate is valid for homogeneous polynomial functionals $J(f) = P(a_n/n)$ with two critical points $z_0, \  z_1$ on $\mathbb S^1$ satisfying $|P(z_0)| > |P(z_1)|$.

\bigskip\noindent
{\bf 4. Generalization}. The above theorems are extended straightforwardly to a more broad class 
$$
\wh S = \bigcup_{|\theta| < \pi, \ |\beta| =1} \wh S_\theta(\beta), 
$$
where $\wh S_\theta(\beta)$ consists of univalent functions
$f(z) = e^{i \theta} z + a_2 z^2 + \dots$ in $\D$ with 
quasiconformal extension to $\hC$ that fix a  point $\beta \in \mathbb S^1$ and their sequential limits under locally uniform convergence in $\D$. 

In this class, the extremal functions  are obtained by post rotations of the Koebe function (4).

\bigskip
{\small\em{ \leftline{Department of Mathematics, Bar-Ilan
University, 5290002 Ramat-Gan, Israel} \leftline{and
Department of Mathematics, University of Virginia,  Charlottesville, VA 22904-4137, USA}}

\end{document}